# Markov Chain properties in terms of column sums of the transition matrix

Jeffrey J. Hunter

*School of Computing and Mathematical Sciences, Auckland University of Technology, Private Bag 92006, Auckland 1142, New Zealand*

September 2, 2011

## Abstract

Questions are posed regarding the influence that the column sums of the transition probabilities of a stochastic matrix (with row sums all one) have on the stationary distribution, the mean first passage times and the Kemeny constant of the associated irreducible discrete time Markov chain. Some new relationships, including some inequalities, and partial answers to the questions, are given using a special generalized matrix inverse that has not previously been considered in the literature on Markov chains.

## 1. Introduction

Let $P = [p_{ij}]$ be the transition matrix of a finite irreducible, discrete time Markov chain $\{X_n\}$, $(n \geq 0)$, with state space $S = \{1, 2,..., m\}$. Such chains have a unique stationary distribution $\{\pi_j\}$, $(1 \leq j \leq m)$ and finite mean first passage times $\{m_{ij}\}$, $(1 \leq i \leq m, 1 \leq j \leq m)$. ([3], [4]).

The stochastic nature of $P$ implies that $\sum_{j=1}^{m} p_{ij} = 1$, for all $i = 1, 2, .., m$, i.e. the row sums are all one.

Let $c_j = \sum_{i=1}^{m} p_{ij}$, for all $j = 1, ..., m$, be the respective column sums of the transition matrix.

We pose the following questions: What influence does the sequence $\{c_j\}$ have on the stationary distribution $\{\pi_j\}$? What influence does the sequence $\{c_j\}$ have on the mean first passage times $\{m_{ij}\}$? Are there relationships connecting the $\{c_j\}$, the $\{\pi_j\}$ and the $\{m_{ij}\}$? Can we deduce bounds on the $\{\pi_j\}$ and the $\{m_{ij}\}$ involving $\{c_j\}$? What effect does the $\{c_j\}$ have on Kemeny's constant $K = \sum_{j=1}^{m} \pi_j m_{ij}$, which is in fact independent of $i \in \{1,2,....,m\}$.

We explore these questions by utilising the generalized matrix inverse $H \equiv [I - P + ec^T]^{-1}$, where $e$ is the column vector of ones and $c^T = (c_1, c_2, ... c_m)$ is the row vector of column sums of $P$. We show that this matrix can also be expressed in terms of Kemeny and Snell's fundamental matrix, or Meyer's group generalized inverse. It has not hitherto been considered as an appropriate matrix to explore the properties of Markov chains.

Two papers appear in the open literature considering the impact that column sums have on the properties of Markov chains. In [11], Kirkland considers the subdominant eigenvalue $\lambda_2$ associated with the set S(c) of n × n stochastic matrices with column sum vector $c^T$. The quantity $\overline{\lambda(c)} = \max\{|\lambda_2(A)| \,|\, A \in S(c)\}$ is considered. The vectors $c^T$ such that $\overline{\lambda(c)} < 1$ are identified and in those cases, nontrivial upper bounds on $\overline{\lambda(c)}$ and weak ergodicity results for forward products are provided. In [12], Kirkland considers an irreducible stochastic matrix $P$ and studies the extent to which the column sum vector for $P$ provides information on a certain condition number $\kappa(P)$, which measures the sensitivity of the stationary distribution vector to perturbations in $P$.

We do not consider these problems in this paper focussing on the main on the key properties of stationary distributions and mean first passage times.

## 2. Properties of the generalized inverse, *H*

In [2] (Theorem 3.3) it was shown that if $\pi^T t \neq 0$ and $u^T e \neq 0$ then $I - P + tu^T$ is non-singular and $[I - P + tu^T]^{-1}$ is a generalized inverse of $I - P$.

For *H*, as defined above, since $\pi^T e = 1$ and $c^T e = e^T P e = e^T e = m$, $I - P + ec^T$ is non-singular and *H* is a generalized inverse of $I - P$.

**Theorem 1**: *If* $H \equiv [I - P + ec^T]^{-1} = [h_{ij}]$,

*then* $c^T H = \pi^T$ *implying* $\pi_j = \sum_{i=1}^{m} c_i h_{ij}$ *for all j = 1, 2, ..., m,* (1)

*and* $He = e/m$ *implying* $h_{i.} \equiv \sum_{j=1}^{m} h_{ij} = 1/m$, *for all i = 1, 2, ..., m,* (2)

*and hence that* $c^T He = 1$, *with* $c^T e = \sum_{i=1}^{m} c_i = m$ . (3)

**Proof:** Since $(I - P + ec^T)H = I$,

$$H - PH + ec^T H = I.$$ (4)

Pre-multiply (4) by $\pi^T$, noting that $\pi^T P = \pi^T$, leads immediately to (1).
Further, since $H(I - P + ec^T) = I$,
$$H - HP + Hec^T = I. \tag{5}$$
Post-multiply (5) by $e$, noting that $Pe = e$, and $c^T e = m$, leads immediately to (2).
Equation (3) follows from (1) and (2). □

The impact of Theorem 1 is that from (1), the stationary probabilities $\pi_j$, for each $j$, can be expressed as a linear function of the column sums of $P$, the $\{c_i\}$, and the elements of $H = [h_{ij}]$.

Alternative simple methods for calculating the stationary probabilities are also given in [7]. The particular method, leading to (1) above, was not included in [7], although a variety of specific forms of the generalized inverses of $I - P$ were considered, leading to various simple expressions for the stationary probabilities.

In order to obtain bounds on the stationary probabilities, we need to have more knowledge of the elements of $H$, in particular the sign of $h_{ij}$. From (2) the sum of the elements in each row of $H$ is $1/m$, but we need more refined information to obtain anything useful.

Let $e_i^T$ be the $i$-th elementary row vector and $e_j$ be the $j$-th elementary column vector.
Let $h_j^{(c)} \equiv He_j$ denote the $j$-th column of $H$ and $h_i^{(r)T} \equiv e_i^T H$ denote the $i$-th row of $H$.
Let $h_{rowsum} = He = \sum_{j=1}^{m} h_j^{(c)} = [h_{1\bullet}, h_{2\bullet}, ..., h_{m\bullet}]^T$ denote the column vector of row sums of $H$ and
$h_{colsum}^T = e^T H = \sum_{j=1}^{m} h_j^{(r)T} = [h_{\bullet 1}, h_{\bullet 2}, ..., h_{\bullet m}]$ denote the row vector of column sums of $H$.

First note that substitution of (1) into (4) yields
$$H - PH + e\pi^T = I. \tag{6}$$
Further substitution of (2) into (5) yields
$$H - HP + ec^T/m = I. \tag{7}$$
Relationships between the rows, columns and elements of $H$ follow from (6) and (7) by pre- and post-multiplication by $e_i^T$ and $e_j$ and the facts that $h_{i\bullet} = h_i^{(r)T} e$, $h_{\bullet j} = e^T h_j^{(c)}$, $h_{ij} = e_i^T He_j$.

The main properties of $H$, deduced from (6) and (7), are summarised in the following theorem.

**Theorem 2**: *The g-inverse $H = [h_{ij}] = [I - P + ec^T]^{-1}$ satisfies the following properties:*

*(a) (Row properties)* $\quad h_i^{(r)T} - p_i^{(r)T} H = e_i^T - \pi^T,$

$\quad\quad\quad\quad\quad\quad\quad\quad\quad\quad h_i^{(r)T} - h_i^{(r)T} P = e_i^T - c^T/m,$

and $\quad\quad\quad\quad\quad\quad\quad\quad h_{i\bullet} = 1/m$.

*(b) (Column properties)* $\quad h_j^{(c)} - Ph_j^{(c)} = e_j - \pi_j e,$

$\quad\quad\quad\quad\quad\quad\quad\quad\quad\quad h_j^{(c)} - Hp_j^{(c)} = e_j - (c_j/m)e,$

$$\text{and} \quad h_{\bullet j} = 1 - (m-1)\pi_j.$$

*(c) (Element properties)*
$$h_{ij} = \sum_{k=1}^{m} p_{ik} h_{kj} + \delta_{ij} - \pi_j,$$
$$h_{ij} = \sum_{k=1}^{m} h_{ik} p_{kj} + \delta_{ij} - c_j/m.$$

*(d) (Row and column sum vectors)*
$$\mathbf{h}_{rowsum} = \mathbf{e}/m,$$
$$\mathbf{h}_{colsum}^T = \mathbf{e}^T - (m-1)\pi^T. \qquad \square$$

While overall row and column sums of the elements of *H* can be derived (Theorem 2(d)), explicit expressions for individual elements of *H* are not readily available.

## 3. Stationary distributions

Irreducible Markov chains have a unique stationary distribution $\{\pi_j\}$, $(1 \leq j \leq m)$. Let $\pi^T = (\pi_1, \pi_2, \ldots, \pi_m)$ be the stationary probability vector of the Markov chain.
The stationary distribution $\{\pi_j\}$, $(j \in S)$ where $S = \{1, 2, \ldots, m\}$ which satisfies the stationary equations:
$$\pi_j = \sum_{i=1}^{m} \pi_i p_{ij} \ (j \in S) \text{ with } \sum_{i=1}^{m} \pi_i = 1. \qquad (8)$$
While we have not solved equations (8) directly, we have shown that $\pi_j = \sum_{i=1}^{m} c_i h_{ij}$ with $\sum_{i=1}^{m} c_i = m$.
Special results for the stationary distribution, in terms of the column sums, are well known in the case of doubly stochastic matrices, (see [9]).

**Theorem 3**: *$\mathbf{c} = \mathbf{e}$ if and only if $\pi = \mathbf{e}/m$.*
Equivalently, $c_1 = c_2 = \ldots = c_m = 1$ *if and only if* $\pi_1 = \pi_2 = \ldots = \pi_m = 1/m$.

**Proof**: From (1), if $\mathbf{c} = \mathbf{e}$, $\pi^T = \mathbf{e}^T H$ that, from Theorem 2(d), implies $\pi^T = \mathbf{e}^T - (m-1)\pi^T$ and hence that $m\pi^T = \mathbf{e}^T$ and the "if" implication follows.
Similarly, if $\pi = \mathbf{e}/m$ then, from (1), $\mathbf{c}^T H = \mathbf{e}^T/m$. Further, from Theorem 2(d), $\mathbf{e}^T H = \mathbf{e}^T - (m-1)\mathbf{e}^T/m = \mathbf{e}^T/m$. Thus $\mathbf{c}^T H = \mathbf{e}^T H$. Post multiplication by $H^{-1} (= I - P + \mathbf{e}\mathbf{c}^T)$, implies $\mathbf{c}^T = \mathbf{e}^T$, leading to the "only if" implication. $\qquad \square$

In order to obtain relationships between the stationary probabilities we need to have information regarding relationships between the elements of *H*.

# 4 Relationship between *H* and *Z*, the fundamental matrix of ergodic Markov chains

In [2] (Theorem 6.3) it was shown that for all generalized inverses *G*, of *I* – *P* that $(I - \Pi)G(I - \Pi)$ is invariant and equals $(I - P)^{\#} = Z - \Pi$, Meyer's group inverse of *I* – *P*, ([13]), where *Z* is Kemeny and Snell's fundamental matrix of ergodic Markov chains, $Z = [I - P + \Pi]^{-1}$, ([10]). We now find an expression for *Z* in terms of *H* and similarly an expression for *H* in terms of *Z*.

**Theorem 4**: If $H = [I - P + ec^T]^{-1}$ and $Z \equiv [I - P + e\pi^T]^{-1}$, then

(a) $\quad Z = H + \Pi - \Pi H$, $\hfill(9)$

(b) $\quad H = Z + \dfrac{1}{m}\Pi - \dfrac{1}{m}ec^T Z$, $\hfill(10)$

(c) $\quad (1+m)\Pi = m\Pi H + ec^T Z$, $\hfill(11)$

(d) $\quad (1+m)\pi^T = m\pi^T H + c^T Z$. $\hfill(12)$

**Proof**: (a) Since $(I - \Pi)H(I - \Pi) = Z - \Pi$, and $H\Pi = He\pi^T = \dfrac{1}{m}e\pi^T = \dfrac{1}{m}\Pi$ simplification yields $Z = (I - \Pi)H(I - \Pi) + \Pi = (I - \Pi)\left(H - \dfrac{1}{m}\Pi\right) + \Pi = (I - \Pi)H + \Pi$, since $(I - \Pi)\Pi = \Pi - \Pi^2 = e\pi^T - e\pi^T e\pi^T = e\pi^T - e\pi^T = 0$, leading to (9).

(d) Premultiplication of (9) by $c^T$ yields $c^T Z = c^T H + c^T e\pi^T - c^T e\pi^T H = \pi^T + m\pi^T - m\pi^T H$ leads to (12).

(c) Premultiplication of (12) by *e* yields (11).

(b) From (11), $\Pi H = \dfrac{1}{m}\left((1+m)\Pi - ec^T Z\right)$, Substitute for $\Pi H$ into (9) yields

$Z = H + \Pi - \left(\dfrac{1}{m}\Pi + \Pi - \dfrac{1}{m}ec^T Z\right) = H - \dfrac{1}{m}\Pi + \dfrac{1}{m}ec^T Z$ leading to (10), an expression for *H* in terms of *Z*. $\hfill\square$

(An alternative proof of Theorem 4 can be given based upon the Sherman-Morrison formula.) Expressions for the elements of *Z* in terms of the elements of *H* and vice-versa can be derived from Theorem 4 and are given in Theorem 5 below. The reason we are interested in these interrelationships is that we know that for ergodic Markov chains the diagonal elements of *Z*, $z_{jj}$, are positive (see below). Matlab examples show that a similar relationship also holds for the diagonal elements of *H*, $h_{jj}$. We seek a theoretical justification for this although a formal proof will need to wait until the next section.

**Theorem 5**: *If* $H = \left[h_{ij}\right]$ *and* $Z = \left[z_{ij}\right]$ *then*

(a) $\quad z_{ij} = h_{ij} + \pi_j - \sum_{k=1}^{m}\pi_k h_{kj}$, $\hfill(13)$

$$(b) \quad h_{ij} = z_{ij} + \frac{1}{m}\pi_j - \frac{1}{m}\sum_{k=1}^{m} c_k z_{kj}, \tag{14}$$

$$(c) \quad (1+m)\pi_j = m\sum_{k=1}^{m}\pi_k h_{kj} + \sum_{k=1}^{m} c_k z_{kj}. \tag{15}$$

□

Note that from (13) or (14), $z_{ij} - h_{ij} = z_{jj} - h_{jj}$, i.e. independent of $i$, so that $z_{ij} - z_{jj} = h_{ij} - h_{jj}$, (See also (20) below for an alternative justification).

The common difference $z_{ij} - h_{ij}$ can be expressed, using (13), as $\pi_j - \sum_{k=1}^{m}\pi_k h_{kj}$ or, by using (14), as $\left(\sum_{k=1}^{m} c_k z_{kj} - \pi_j\right)/m$ leading also to (15), which is the elemental expression of (12).

In [5], parametric forms of generalized inverses of Markovian kernels were derived with every one condition g-inverse $G$ of $I - P$ expressed in the parametric form $G(\alpha,\beta,\gamma) = [I - P + \alpha\beta^T]^{-1} + \gamma e\pi^T$, where $\alpha = [I - (I - P)G]e$, $\beta^T = \pi^T[I - G(I - P)]$ and $\gamma + 1 = \beta^T G\alpha$.

$Z = [I - P - e\pi^T]^{-1}$ is actually expressed in parametric form with $\alpha = e$, $\beta^T = \pi^T$ and $\gamma = 0$.

For $H$ it is easily verified that $H = [I - P + ec^T]^{-1} = \left[I - P + \frac{1}{m}ee^T\right]^{-1} + \left(\frac{1}{m} - 1\right)e\pi^T$. Thus, if necessary, $H$ can be re-expressed in a form involving $\pi^T$, as well as in terms of $Z$, as above.

**Theorem 6**: If $c_k = 1$ for all $k$,

$$H = Z + \left(\frac{1-m}{m^2}\right)E. \tag{16}$$

**Proof**: From Theorem 3 $\pi_k = 1/m$ for all $k$, and hence from Theorem 2(b) and Theorem 5 (c), $z_{gj} = 1, h_{gj} = \frac{1}{m}$. Thus, from Theorem 5(b,) $h_{ij} = z_{ij} + \frac{1-m}{m^2}$ leading to expression (16). □

## 5. Mean first passage times

Let $M = [m_{ij}]$ be the matrix of expected first passage times from state $i$ to state $j$ in an irreducible finite Markov chain with transition matrix $P$. The following result is well known. (See for example, [2] (Section 5.1), [4] (Corollary 7.3.3B), [10] (Theorem 4.4.4).)

**Theorem 7**: *M satisfies the matrix equation*

$$(I - P)M = E - PM_d, \tag{17}$$

*where $E = ee^T = [1]$, $M_d = [\delta_{ij}m_{ij}]$, a diagonal matrix with elements the diagonal elements of M. Further, $M_d = (\Pi_d)^{-1} \equiv D$, where $\Pi = e\pi^T$.* □

It is well known that the solution of equations of the form of (16) can be effected using g-inverses of $I - P$, (see e.g. [2] and [4]). Any g-inverse of $I - P$ has the form $G = \left[I - P + tu^T\right]^{-1} + ef^T + g\pi^T$, where $u^T e \neq 0$, $\pi^T t \neq 0$ and $f$ and $g$ are arbitrary vectors. The following result (17) appears in [2] (Theorem 5.1) and [4] (Theorem 7.3.6) and while result (18)) appears in [8] (Corollary 2.3.2).

**Theorem 8**: *If G is any g-inverse of I – P, then*

$$M = [G\Pi - E(G\Pi)_d + I - G + EG_d]D. \tag{18}$$

*Further, under any of the following three equivalent conditions,*
   (i)     $Ge = ge$, $g$ a constant,
   (ii)    $GE - E(G\Pi)_d D = 0$,
   (iii)   $G\Pi - E(G\Pi)_d = 0$,

$$M = [I - G + EG_d]D. \tag{19}$$

□

As a result of (2), $H$ as defined satisfies condition (*i*) of Theorem 8 so that expression (19) is valid for $G = H$. Other special cases for equation (18)) are $G = Z$, Kemeny and Snell's fundamental matrix $Z = [I - P - \Pi]^{-1}$ (since $Ze = e$ and $g = 1$, as given initially in [10] (Theorem 4.4.7)) and $G = = (I - P)^\# = Z - \Pi$, Meyer's group inverse of $I - P$, (with $(I - P)^\# e = 0$ and $g = 0$) as given in [13] (Theorem 3.3).

Elemental expressions for the $m_{ij}$ follow from Theorem 8, as follows:

**Theorem 9**: *If $G = [g_{ij}]$ is any generalized inverse of $I - P$, then*
$$m_{ij} = ([g_{jj} - g_{ij} + \delta_{ij}]/\pi_j) + (g_{i\bullet} - g_{j\bullet}), \text{ for all } i, j. \tag{20}$$
*Further, when $Ge = ge$,*
$$m_{ij} = [g_{jj} - g_{ij} + \delta_{ij}]/\pi_j, \text{ for all } i, j. \tag{21}$$

**Proof**: Expressing (18) and (19) in elemental form leads to (20) and (21), respectively. □

We have some key observations from Theorem 9. Since $H = [h_{ij}]$ satisfies the required condition for (21), we have that

$$m_{ij} = \begin{cases} \dfrac{1}{\pi_j} = \dfrac{1}{\sum_{i=1}^m c_i h_{ij}}, & i = j, \\ \dfrac{h_{jj} - h_{ij}}{\pi_j} = \dfrac{h_{jj} - h_{ij}}{\sum_{i=1}^m c_i h_{ij}}, & i \neq j. \end{cases} \tag{22}$$

Thus a knowledge of the column sums $\{c_i\}$ and the elements $h_{ij}$ leads directly to expressions for the mean first passage times.

Relationships between the $\{\pi_j\}$, $\{m_{ij}\}$ and $\{c_j\}$ can be derived from Theorems 7 and 8.

**Theorem 10**: *For all $j \in \{1, 2, ..., m\}$*,

$$\sum_{i=1}^{m} m_{ij} - \sum_{i=1}^{m} c_i m_{ij} = m - \frac{c_j}{\pi_j} = m - c_j m_{jj}, \quad (23)$$

$$\sum_{i=1}^{m} c_i m_{ij} = \frac{c_j}{\pi_j} - 1 + \frac{mh_{jj}}{\pi_j} = c_j m_{jj} - 1 + mh_{jj} m_{jj} \quad (24)$$

*and hence*

$$\sum_{i=1}^{m} m_{ij} = m_{\bullet j} = m - 1 + \frac{mh_{jj}}{\pi_j} = m - 1 + mh_{jj} m_{jj}, \quad (25)$$

**Proof**: Pre-multiplication of (17) by $e^T$, since $e^T P = c^T$, yields $e^T M - c^T M = me^T - c^T M_d$. Further, $e^T M = (m_{\bullet 1}, ..., m_{\bullet j}, ..., m_{\bullet m})$, where $\boxed{m_{\bullet j} = \sum_{i=1}^{m} m_{ij}}$, and $c^T M_d = \left(\frac{c_1}{\pi_1}, ..., \frac{c_j}{\pi_j}, ..., \frac{c_m}{\pi_m}\right)$, and (23) follows upon extracting the *j*-th element.

From (18) with $G = H$, pre-multiplication by $c^T$ gives $c^T M = c^T D - c^T HD + c^T ee^T H_d D$. Observing that $c^T M = (\sum_{i=1}^{m} c_i m_{i1}, ..., \sum_{i=1}^{m} c_i m_{ij}, ..., \sum_{i=1}^{m} c_i m_{im})$, that $c^T D = c^T M_d$ as above, $c^T HD = \pi^T M_d = (1, ..., 1, ..., 1) = e^T$, and $c^T ee^T H_d D = me^T H_d M_d = m\left(\frac{h_{11}}{\pi_1}, ..., \frac{h_{jj}}{\pi_j}, ..., \frac{h_{mm}}{\pi_m}\right)$, it is easily seen, upon extracting the *j*-th element in the combined expression, that (24) follows.

Expression (25) follows from (23) upon substitution of $\sum_{i=1}^{m} c_i m_{ij}$, from (24). □

Expression (23) gives a new interesting connection between the $\{c_j\}$ and the $\{m_{ij}\}$, unhindered by the elements of *H*. A simple extension of Theorem 10, re-expressing (23) and (24) in terms of the stationary probability $\pi_j$, yields the following new results, where the second expressions of (26), (27) and (28) follow by multiplying out the first expressions and using the observation that $\pi_j m_{jj} = 1$,

$$\pi_j = \frac{c_j}{m - \sum_{i=1}^{m} m_{ij} + \sum_{i=1}^{m} c_i m_{ij}} = \frac{1}{m - \sum_{i \neq j} m_{ij} + \sum_{i \neq j} c_i m_{ij}}, \quad (26)$$

$$\pi_j = \frac{c_j + mh_{jj}}{1 + \sum_{i=1}^{m} c_i m_{ij}} = \frac{mh_{jj}}{1 + \sum_{i \neq j} c_i m_{ij}}, \quad (27)$$

$$\pi_j = \frac{mh_{jj}}{1 + \sum_{i=1}^{m} m_{ij} - m} = \frac{mh_{jj} - 1}{1 + \sum_{i \neq j} m_{ij} - m}. \quad (28)$$

Now note that from (22), and (21), that since *Z*, as well as *H*, also satisfies the required condition for (21), with $Ze = e$, we have that $\pi_j m_{ij} = h_{jj} - h_{ij} = z_{jj} - z_{ij} > 0$ for all $i \neq j$.

From [4] (Theorem 7.3.8), $\pi^T M = e^T Z_d D = (z_{11}/\pi_1, ..., z_{jj}/\pi_j, ..., z_{mm}/\pi_m)$ so that it is clear that since all the terms on the left hand side are positive we deduce that $z_{jj} > 0$ for all $j$.

**Theorem 11:** *For all $i, j \in \{1, 2, ..., m\}$,*

$$h_{ij} = \begin{cases} \pi_j (1 + \sum_{k \neq j} c_k m_{kj})/m, & i = j, \\ \pi_j (1 - m m_{ij} + \sum_{k \neq j} c_k m_{kj})/m, & i \neq j. \end{cases} \quad (29)$$

implying $\quad H = \dfrac{1}{m} \Pi + (\dfrac{1}{m} C - I)(M - M_d)\Pi_d \quad$ where $\Pi = e\pi^T$ and $C = ec^T$.

**Proof:** The expression for $h_{jj}$ follows, from (27) while the expression for $h_{ij}$ follows from the observation that $h_{ij} = h_{jj} - \pi_j m_{ij}$ for $i \neq j$. The expression for $H$ follows from (29). □

As a further observation note that $CH = \Pi, C\Pi = m\Pi$ and $C^2 = mC$.

An important consequence of (29) is that $h_{jj} > 0$. The usefulness of this is that we now know that in the key relationship $\pi_j = \sum_{i=1}^m c_i h_{ij}$, that the coefficient of $c_j$ is always positive. However we have no surety regarding the positivity of any of the remaining terms $h_{ij}$ when $i \neq j$. We explore some further consequences in the next section.

From (26) and (27) we have some new bounds for the stationary probability $\pi_j$

$$\pi_j > max\left[\frac{1}{m + \sum_{i \neq j} c_i m_{ij}}, \frac{c_j}{1 + \sum_{i=1}^m c_i m_{ij}}\right].$$

For all finite irreducible Markov chains $\pi_i \leq \pi_j \Leftrightarrow m_{jj} \leq m_{ii}$ since $m_{ii} = 1/\pi_i$.

**Theorem 12:** If $c_i = 1$ for all $i$, then
$$m_{\bullet j} = m - 1 + m^2 h_{jj} = m^2 z_{jj}. \quad (30)$$

**Proof:** From Theorem 3, $\pi_i = 1/m$ for all $i$, and the first expression of (30) follows from (25). From (16), $m^2 h_{jj} = m^2 z_{jj} + 1 - m$, so that $m_{\bullet j} = m^2 z_{jj}$. □

This is an extension of the results derived in [9] for doubly stochastic Markov chains. We explore results for $m_{ig}$ in the next section.

Thus for the situation of constant column sums, $m_{\bullet j} \leq m_{\bullet i} \Leftrightarrow z_{jj} \leq z_{ii} \Leftrightarrow h_{jj} \leq h_{ii}$

## 6. Kemeny's constant

Kemeny [10] made the observation that the following expression $K_i = \sum_{j=1}^{m} m_{ij}\pi_j$ is in fact independent of $i \in \{1,2,....,m\}$, so that $K_i = K$. It has since been realised that this constant has many important interpretations in terms of properties of the Markov chain, in particular being used in the properties of mixing (expected time to stationarity) and a constant used in bounding overall differences in the stationary probabilities of a chain subjected to perturbations ([6]).

We have the following general results, initially derived in [6] (Theorem 2.4) (See also [8], Theorem 3.2,) that if $G = [g_{ij}]$ is any g-inverse of $I - P$, then Kemeny's constant $K$ has the form

$$K = 1 + tr(G) - tr(G\Pi) = 1 + \sum_{j=1}^{m}(g_{jj} - g_{j\bullet}\pi_j). \tag{31}$$

This leads to the following specific expressions for $K$:

**Theorem 13**:
$$K = 1 - (1/m) + tr(H) = 1 - (1/m) + \sum_{j=1}^{m} h_{jj}, \tag{32}$$

$$= tr(Z) = \sum_{j=1}^{m} z_{jj}. \tag{33}$$

**Proof**: If for some $g$, $Ge = ge$, then from (31) above $K = 1 - g + tr(G)$.
$G = H$ and $G = Z$ both have the required property of $Ge = ge$, since $He = 1/m \, e$ and $Ze = e$, and (32) and (33) both follow from (31) with $g = 1/m$ and 1, respectively. □

Substitution for $h_{jj}$ from the second expression of (27) into (32) leads to

$$K = 1 - \frac{1}{m} + \frac{1}{m}\sum_{j=1}^{m}\pi_j(1 + \sum_{k\neq j}c_k m_{kj}),$$

$$= 1 - \frac{1}{m} + \frac{1}{m}\sum_{j=1}^{m}\pi_j + \frac{1}{m}\sum_{j=1}^{m}\pi_j(\sum_k c_k m_{kj}) - \frac{1}{m}\sum_{j=1}^{m}c_j,$$

$$= \frac{1}{m}\sum_{j=1}^{m}\sum_{k=1}^{m}c_k\pi_j m_{kj} = \frac{1}{m}\sum_{k=1}^{m}c_k\sum_{j=1}^{m}\pi_j m_{kj}.$$

Note however that this follows directly from the definition of $K_i = \sum_{j=1}^{m}\pi_j m_{ij}$ and the facts that $K_i = K$ and $\sum_{k=1}^{m}c_k = m$.

In [4] it was shown that for any irreducible m-state Markov chain that $K \geq \frac{m+1}{2}$, so that from (32), after simplification,

$$tr(H) = \sum_{j=1}^{m} h_{jj} \geq \frac{m-1}{2} + \frac{1}{m}. \tag{34}$$

Further we have the following key properties: $h_{jj} > 0$ (from (27)), and, for $i \neq j$, $h_{jj} > h_{ij}$ (from (22)), with $\sum_{j=1}^{m} h_{ij} = 1/m$, (from (2)). These results lead to the following observation.

**Theorem 14**: For all $j \in S = \{1,2,....,m\}$,
$$\pi_j < mh_{jj}. \tag{35}$$

**Proof**: From (1) $\pi_j = \sum_{i=1}^{m} c_i h_{ij} = c_1 h_{1j} + ... + c_j h_{jj} + ... + c_m h_{mj}$

$< c_1 h_{jj} + ... + c_j h_{jj} + ... + c_m h_{jj} = (\sum_{k=1}^{m} c_k) h_{jj} = mh_{jj}$ using (22) and (3).

This result also follows from (27) since $h_{jj} = \pi_j (1 + \sum_{k \neq j} c_k m_{kj})/m > \pi_j/m$. □

Result (35) implies that $m \sum_{j=1}^{m} h_{jj} > \sum_{j=1}^{m} \pi_j = 1$ and hence that $\sum_{j=1}^{m} h_{jj} > 1/m$, a much weaker bound that that given by (34).

**Theorem 15:** If $c_i = 1$ for all $i$, then
$$m_{i\bullet} = mK = m - 1 + mtr(H) = mtr(Z). \tag{36}$$

Further
$$K = \frac{m_{\bullet\bullet}}{m^2} = tr(Z). \tag{37}$$

**Proof**: The conditions imply a uniform stationary distribution so that $K_i = \sum_{j=1}^{m} m_{ij} \pi_j = \frac{m_{i\bullet}}{m}$
$= K$ for all $i$, so that the mean first passage time matrix has constant row sums given by (36), as also observed in [9].

Further $m_{\bullet\bullet} = mm_{i\bullet} = m^2 K$. In this case, from (30), $m_{\bullet j} = m^2 z_{jj}$ so that $m_{\bullet\bullet} = \sum_{j=1}^{m} m_{\bullet j} =$
$m^2 \sum_{j=1}^{m} z_{jj} = m^2 tr(Z)$ leading to (37), (consistent with (33)). □

Note also for constant column sums of the transition matrix that, from (36) and the lower bound on $K$, that for all $i$, $m_{i\bullet} \geq \frac{m(m+1)}{2}$, a new result.

We illustrate some of these properties with a series of examples.

## 7. Examples

**Example 1: (Two-state Markov chain)**

Let $P = \begin{bmatrix} p_{11} & p_{12} \\ p_{21} & p_{22} \end{bmatrix} = \begin{bmatrix} 1-a & a \\ b & 1-b \end{bmatrix}$, $(0 \leq a \leq 1, 0 \leq b \leq 1)$, be the transition matrix of a two-state Markov chain with state space $S = \{1, 2\}$. Let $d = 1 - a - b$ so that $1 - d = a + b$.

To ensure that the Markov chain is irreducible we henceforth assume that $-1 \leq d < 1$ so that the Markov chain has a unique stationary probability vector given by, ([4], p.71),
$$\pi^T = (\pi_1, \pi_2) = (b/(a+b), a/(a+b)).$$
(If $-1 < d < 1$, the Markov chain is regular and the stationary distribution is in fact the limiting distribution. If $d = -1$ the Markov chain is irreducible periodic, period 2.)

The row vector of column sums is given by $c^T = (c_1, c_2) = (1-(a-b), 1+(a-b))$

Note that where as the parameters $a$ and $b$ specify all the transition probabilities, the parameters $c_1$ and $c_2$ do not uniquely specify the transition probabilities since $c_1 + c_2 = 2$ with $c_2 - c_1 = 2(a-b)$, and we cannot solve for $a$ and $b$ in terms of $c_1$ and $c_2$.

The matrix of mean first passage times is given by
$$M = \begin{bmatrix} m_{11} & m_{12} \\ m_{21} & m_{22} \end{bmatrix} = \begin{bmatrix} (a+b)/b & 1/a \\ 1/b & (a+b)/a \end{bmatrix}, \text{([4], p. 135), or ([10], p. 94))},$$

The $H$ matrix is given by $H = [I - P + ec^T]^{-1} = \dfrac{1}{2(a+b)} \begin{bmatrix} 1+a & -(1-b) \\ -(1-a) & 1+b \end{bmatrix}$.

The fundamental matrix $Z$ is given by, ([4], p. 135), $Z = \dfrac{1}{a+b} \begin{bmatrix} b + \dfrac{a}{a+b} & a - \dfrac{a}{a+b} \\ b - \dfrac{b}{a+b} & a + \dfrac{b}{a+b} \end{bmatrix}$.

Kemeny's constant is given by, ([6]), $K = 1 + \dfrac{1}{a+b}$, with the property that $K \geq 1.5$ with the minimum value occurring when $a = b = 1$.

It is easy to verify that for this two-state case
$$c_1 \leq c_2 \Leftrightarrow b \leq a \Leftrightarrow \pi_1 \leq \pi_2 \Leftrightarrow m_{22} \leq m_{11},$$
and
$$h_{11} \leq h_{22} \Leftrightarrow a \leq b \Leftrightarrow m_{11} + m_{21} = m_{\cdot 1} \leq m_{\cdot 2} = m_{12} + m_{22}.$$

This suggests some possibilities that need to be explored for larger state spaces.

**Example 2: (Three-state Markov chain)**

Let $P = \begin{bmatrix} p_{11} & p_{12} & p_{13} \\ p_{21} & p_{22} & p_{23} \\ p_{31} & p_{32} & p_{33} \end{bmatrix} = \begin{bmatrix} 1 - p_2 - p_3 & p_2 & p_3 \\ q_1 & 1 - q_1 - q_3 & q_3 \\ r_1 & r_2 & 1 - r_1 - r_2 \end{bmatrix}$ be the transition matrix of a Markov chain with state space $S = \{1, 2, 3\}$, where we have used the parametrisation as used in [6], [8]. Note that for the six constrained parameters we have $0 < p_2 + p_3 \leq 1$, $0 < q_1 + q_3 \leq 1$ and $0 < r_1 + r_2 \leq 1$.

Let $\Delta_1 \equiv q_3 r_1 + q_1 r_2 + q_1 r_1$, $\Delta_2 \equiv r_1 p_2 + r_2 p_3 + r_2 p_2$, $\Delta_3 \equiv p_2 q_3 + p_3 q_1 + p_3 q_3$, $\Delta \equiv \Delta_1 + \Delta_2 + \Delta_3$.

The Markov chain, with the above transition matrix, is irreducible (and hence a stationary distribution exists) if and only if $\Delta_1 > 0$, $\Delta_2 > 0$, $\Delta_3 > 0$.

It is easily shown that the stationary probability vector is given by $(\pi_1, \pi_2, \pi_3) = \frac{1}{\Delta}(\Delta_1, \Delta_2, \Delta_3)$.

Define $\tau_{12} = p_3 + r_1 + r_2$, $\tau_{13} = p_2 + q_1 + q_3$, $\tau_{21} = q_3 + r_1 + r_2$, $\tau_{23} = q_1 + p_2 + p_3$, $\tau_{31} = r_2 + q_1 + q_3$, $\tau_{32} = r_1 + p_2 + p_3$, $\tau_. = p_2 + p_3 + q_1 + q_3 + r_1 + r_2$, so that $\tau_. = \tau_{12} + \tau_{13} = \tau_{21} + \tau_{23} = \tau_{31} + \tau_{32}$.

In [6], a general expression for any generalized inverse of $I - P$ of the form $[I - P + tu^T]^{-1}$ was given. In particular

$$G(e,u) = [I - P + eu^T]^{-1} = \frac{1}{u_.}\left[\Pi + u_1 A_1 + u_2 A_2 + u_3 A_3\right] \text{ where } u_. = u_1 + u_2 + u_3,$$

$$\Pi = \frac{1}{\Delta}\begin{bmatrix} \Delta_1 & \Delta_2 & \Delta_3 \\ \Delta_1 & \Delta_2 & \Delta_3 \\ \Delta_1 & \Delta_2 & \Delta_3 \end{bmatrix}, A_1 = \frac{1}{\Delta}\begin{bmatrix} 0 & 0 & 0 \\ -\tau_{21} & \tau_{12} & \tau_{21} - \tau_{12} \\ -\tau_{31} & \tau_{31} - \tau_{13} & \tau_{13} \end{bmatrix}, A_2 = \frac{1}{\Delta}\begin{bmatrix} \tau_{21} & -\tau_{12} & \tau_{12} - \tau_{21} \\ 0 & 0 & 0 \\ \tau_{32} - \tau_{23} & -\tau_{32} & \tau_{23} \end{bmatrix}$$

and $A_3 = \frac{1}{\Delta}\begin{bmatrix} \tau_{31} & \tau_{13} - \tau_{31} & -\tau_{13} \\ \tau_{23} - \tau_{32} & \tau_{32} & -\tau_{23} \\ 0 & 0 & 0 \end{bmatrix}$.

From the above general expression, since $H = G(e, c)$ and $Z = G(e, \pi)$,

$$H = \frac{1}{3}\left[\Pi + c_1 A_1 + c_2 A_2 + c_3 A_3\right]$$

and

$$Z = \Pi + \pi_1 A_1 + \pi_2 A_2 + \pi_3 A_3.$$

Since $u^T G = \pi^T$ and since $u^T \Pi = u^T e \pi^T = u_. \pi^T$ we must have $u^T(u_1 A_1 + u_2 A_2 + u_3 A_3) = 0^T$. Thus $c^T H = \pi^T$ (taking $u^T = c^T$ with $c^T e = c_g = 3$) and $\pi^T Z = \pi^T$ (taking $u^T = \pi^T$ with $\pi^T e = \pi_. = 1$) as established for all irreducible finite state cases.

It is easily verified, by direct multiplication, that for all $(u_1, u_2, u_3)$,

$$(u_1, u_2, u_3)\begin{bmatrix} u_2 \tau_{21} + u_3 \tau_{31} & -u_2 \tau_{12} + u_3 \tau_{13} - u_3 \tau_{31} & u_2 \tau_{12} - u_2 \tau_{21} - u_3 \tau_{13} \\ -u_1 \tau_{21} + u_3 \tau_{23} - u_3 \tau_{32} & u_1 \tau_{12} + u_3 \tau_{32} & u_1 \tau_{21} - u_1 \tau_{12} - u_3 \tau_{23} \\ -u_1 \tau_{31} + u_2 \tau_{32} - u_2 \tau_{23} & u_1 \tau_{31} - u_1 \tau_{13} - u_2 \tau_{32} & u_1 \tau_{13} u_2 \tau_{23} \end{bmatrix} = (0,0,0).$$

independently establishing that $u^T(u_1 A_1 + u_2 A_2 + u_3 A_3) = 0^T$.

Upon substitution,

$$H = \frac{1}{3\Delta}\begin{bmatrix} \Delta_1 + c_2 \tau_{21} + c_3 \tau_{31} & \Delta_2 - c_2 \tau_{12} + c_3(\tau_{13} - \tau_{31}) & \Delta_3 + c_2(\tau_{12} - \tau_{21}) - c_3 \tau_{13} \\ \Delta_1 - c_1 \tau_{21} + c_3(\tau_{23} - \tau_{32}) & \Delta_2 + c_1 \tau_{12} + c_3 \tau_{32} & \Delta_3 + c_1(\tau_{21} - \tau_{12}) - c_3 \tau_{23} \\ \Delta_1 - c_1 \tau_{31} + c_2(\tau_{32} - \tau_{23}) & \Delta_2 + c_1(\tau_{31} - \tau_{13}) - c_2 \tau_{32} & \Delta_3 + c_1 \tau_{13} + c_2 \tau_{23} \end{bmatrix}$$

and

$$Z = \frac{1}{\Delta} \begin{bmatrix} \Delta_1 + \Delta_2\tau_{21} + \Delta_3\tau_{31} & \Delta_2(1-\tau_{12}) + \Delta_3(\tau_{13}-\tau_{31}) & \Delta_2(\tau_{12}-\tau_{21}) + \Delta_3(1-\tau_{13}) \\ \Delta_1(1-\tau_{21}) + \Delta_3(\tau_{23}-\tau_{32}) & \Delta_1\tau_{12} + \Delta_2 + \Delta_3\tau_{32} & \Delta_1(\tau_{21}-\tau_{12}) + \Delta_3(1-\tau_{23}) \\ \Delta_1(1-\tau_{31}) + \Delta_2(\tau_{32}-\tau_{23}) & \Delta_1(\tau_{31}-\tau_{13}) + \Delta_2(1-\tau_{32}) & \Delta_1\tau_{13} + \Delta_2\tau_{23} + \Delta_3 \end{bmatrix}.$$

Note that it easily seen by examining the diagonal elements on the above two matrices that $h_{ii} > 0$ and $z_{ii} > 0$ for all $i$, as established in general in (28). The non-negativity of $h_{jj} - h_{ij}$ and $z_{jj} - z_{ij}$ for all $i \neq j$ leads (using (22) or the equivalent expression using the $z_{ij}$ to the following expression for the mean first passage time matrix, (see also [6]):

$$M = \begin{bmatrix} \dfrac{\Delta}{\Delta_1} & \dfrac{\tau_{12}}{\Delta_2} & \dfrac{\tau_{13}}{\Delta_3} \\ \dfrac{\tau_{21}}{\Delta_1} & \dfrac{\Delta}{\Delta_2} & \dfrac{\tau_{23}}{\Delta_3} \\ \dfrac{\tau_{31}}{\Delta_1} & \dfrac{\tau_{32}}{\Delta_2} & \dfrac{\Delta}{\Delta_3} \end{bmatrix}.$$

The expected "time to mixing" is given, (see [6]), as $K = 1 + (\tau/\Delta)$.

In [6] it was shown that for all three-state irreducible Markov chains, $K \geq 2$, (with $K = 2$ achieved in "the minimal period 3" case when $p_2 = q_3 = r_1$).

We now explore some possible relationships alluded to in the two-state case.

Under the imposition of column totals with $c_1 + c_2 + c_3 = 3$, we can reduce the free parameters to $p_2$, $p_3$, $q_1$, $q_3$, $c_1$ and $c_2$ by taking $r_1 = c_1 - 1 + p_2 + p_3 - q_1$, $r_2 = c_2 - 1 - p_2 + q_1 + q_3$.

Let $\alpha_1 \equiv q_1 + q_3 - p_2$, $\alpha_2 \equiv p_2 + p_3 - q_1$, then
$\pi_1 \leq \pi_2 \Leftrightarrow m_{22} \leq m_{11} \Leftrightarrow \Delta_1 \leq \Delta_2 \Leftrightarrow q_3 r_1 + q_1 r_2 + q_1 r_1 \leq r_1 p_2 + r_2 p_3 + r_2 p_2 \Leftrightarrow r_1 \alpha_1 \leq r_2 \alpha_2$,
$c_1 \leq c_2 \Leftrightarrow 2q_1 + q_3 + r_1 \leq 2p_2 + p_3 + r_2 \Leftrightarrow r_1 + \alpha_1 \leq r_2 + \alpha_2$.

Unfortunately we cannot deduce universal, if and only if, inequalities connecting $c_1 \leq c_2$ with $\pi_1 \leq \pi_2$. The following table gives parameter regions where the stated inequalities occur, in the case where $r_1 > 0$.

|  | $c_1 \leq c_2$ | $c_1 \geq c_2$ |  |
|---|---|---|---|
| $\pi_1 \leq \pi_2$ | $\alpha_1 \leq \min\left(\left(\dfrac{r_2}{r_1}\right)\alpha_2, r_2 - r_1 + \alpha_2\right)$ | $r_2 - r_1 + \alpha_2 \leq \alpha_1 \leq \left(\dfrac{r_2}{r_1}\right)\alpha_2$ | $r_1\alpha_1 \leq r_2\alpha_2.$ |
| $\pi_1 \geq \pi_2$ | $\left(\dfrac{r_2}{r_1}\right)\alpha_2 \leq \alpha_1 \leq r_2 - r_1 + \alpha_2$ | $\alpha_1 \geq \max\left(\left(\dfrac{r_2}{r_1}\right)\alpha_2, r_2 - r_1 + \alpha_2\right)$ | $r_2\alpha_2 \leq r_1\alpha_1$ |
|  | $r_1 + \alpha_1 \leq r_2 + \alpha_2$ | $r_2 + \alpha_2 \leq r_1 + \alpha_1$ |  |

When $r_1 = 0$ we require for irreducibility that $\Delta_1 = q_1 r_2 > 0$ implying $r_2 > 0$ and $q_1 > 0$. Further, we require $\Delta_2 = (p_2 + p_3)r_2 > 0$ so that $(p_2 + p_3) > 0$. Note also that we require $\Delta_3 \equiv (p_2 + p_3)q_3 + p_3 q_1 > 0$ so that either $q_3 > 0$ and/or $p_3 > 0$. Now $\alpha_1 \equiv q_1 + q_3 - p_2$, $\alpha_2 \equiv p_2 + p_3 - q_1$. Thus $0 \le \alpha_2 \Leftrightarrow q_1 \le p_2 + p_3$ and $\alpha_1 \le r_2 + \alpha_2 \Leftrightarrow 2q_1 + q_3 \le r_2 + 2p_2 + p_3$.

|  | $c_1 \le c_2$ | $c_1 \ge c_2$ |  |
|---|---|---|---|
| $\pi_1 \le \pi_2$ | $\alpha_1 \le r_2 + \alpha_2, 0 \le \alpha_2$ | $r_2 + \alpha_2 \le \alpha_1, 0 \le \alpha_2$ | $0 \le \alpha_2$ |
| $\pi_1 \ge \pi_2$ | $\alpha_1 \le r_2 + \alpha_2, \alpha_2 \le 0$ | $r_2 + \alpha_2 \le \alpha_1, \alpha_2 \le 0$ | $\alpha_2 \le 0$ |
|  | $\alpha_1 \le r_2 + \alpha_2$ | $r_2 + \alpha_2 \le \alpha_1$ |  |

In terms of relationships between the $m_{\cdot 1}$ and $m_{\cdot 2}$ we do not have any inequalities expressed in terms of the column sums, but in terms of the stationary probabilities:

$$m_{\cdot 2} \le m_{\cdot 1} \Leftrightarrow \Delta_1(\Delta + \tau_{12} + \tau_{32}) \le \Delta_2(\Delta + \tau_{21} + \tau_{31})$$
$$\Leftrightarrow \Delta(\Delta_1 - \Delta_2) \le \Delta_2(\tau_{21} + \tau_{31}) - \Delta_1(\tau_{12} + \tau_{32})$$
$$\Leftrightarrow \Delta_1 - \Delta_2 \le \pi_2(\tau_{21} + \tau_{31}) - \pi_1(\tau_{12} + \tau_{32}).$$

It is possible that these inequalities are related to inequalities involving the diagonal elements of the $Z$ and or $H$ matrix, as below, but no obvious universal bounds appear:

$$z_{11} \le z_{22} \Leftrightarrow \Delta_1 + \pi_2\tau_{21} + \pi_3\tau_{31} \le \Delta_2 + \pi_1\tau_{12} + \pi_3\tau_{32}$$
$$\Leftrightarrow \Delta_1 - \Delta_2 \le \pi_1\tau_{12} - \pi_2\tau_{21} + \pi_3(\tau_{32} - \tau_{31})$$
$$\Leftrightarrow \Delta(\Delta_1 - \Delta_2) \le \Delta_1\tau_{12} - \Delta_2\tau_{21} + \Delta_3(\tau_{32} - \tau_{31}).$$

$$h_{11} \le h_{22} \Leftrightarrow \Delta_1 + c_2\tau_{21} + c_3\tau_{31} \le \Delta_2 + c_1\tau_{12} + c_3\tau_{32}$$
$$\Leftrightarrow \Delta_1 - \Delta_2 \le c_1\tau_{12} - c_2\tau_{21} + c_3(\tau_{32} - \tau_{31}).$$

**Example 3: (Five state Markov chain)**

We consider the following irreducible five state Markov chain example taken from Kemeny and Snell [10] (p199) with transition matrix given by

$$\begin{bmatrix} 0.831 & 0.033 & 0.013 & 0.028 & 0.095 \\ 0.046 & 0.788 & 0.016 & 0.038 & 0.112 \\ 0.038 & 0.034 & 0.785 & 0.036 & 0.107 \\ 0.054 & 0.045 & 0.017 & 0.728 & 0.156 \\ 0.082 & 0.065 & 0.023 & 0.071 & 0.759 \end{bmatrix}$$

The column sum vector for this transition matrix is given by $(1.051, 0.965, 0.854, 0.910, 1.229)$. Let us rearrange the states in the order $\{5, 1, 2, 4, 3\}$ so that the column sums are ordered as $c^T = (c_1, c_2, c_3, c_4, c_5)$ with $c_1 \geq c_2 \geq c_3 \geq c_4 \geq c_5$ so that $c^T = (1.229, 1.051, 0.965, 0.910, 0.854)$. The transition matrix for this particular Markov chain is then represented as

$$P = \begin{bmatrix} 0.759 & 0.082 & 0.065 & 0.071 & 0.023 \\ 0.095 & 0.831 & 0.033 & 0.028 & 0.013 \\ 0.112 & 0.046 & 0.788 & 0.038 & 0.016 \\ 0.156 & 0.054 & 0.045 & 0.728 & 0.017 \\ 0.107 & 0.038 & 0.034 & 0.036 & 0.785 \end{bmatrix}.$$

Henceforth, we consider the properties of the chain with this reordered state space relabelled as states $\{1, 2, 3, 4, 5\}$.

The stationary probability vector is given by $\pi^T = (0.3216, 0.2705, 0.1842, 0.1476, 0.0761)$ implying that $\pi_1 \geq \pi_2 \geq \pi_3 \geq \pi_4 \geq \pi_5$. Thus the stationary probabilities appear as in the same order as in the vector of column sums. This is not the result that we necessarily expected.

$$H = \begin{bmatrix} 2.1984 & -0.5537 & -0.4911 & -0.3007 & -0.6530 \\ -0.8883 & 3.5691 & -0.9174 & -0.7613 & -0.8021 \\ -0.6457 & -1.0375 & 3.2047 & -0.5873 & -0.7342 \\ -0.2485 & -0.8505 & -0.6652 & 2.6746 & -0.7104 \\ -0.7023 & -1.2092 & -0.8680 & -0.6157 & 3.5952 \end{bmatrix}.$$

This $H$ matrix has the property that all the diagonal elements are positive (as expected by result (28)). The off-diagonal terms are all negative, although this is not an expected result in general (see the eight state example to follow.) Each row sum is 0.200, consistent with (2) in Theorem 1 and the column sums are given as $h^T_{colsum} = (-0.2863, -0.0818, 0.2631, 0.4096, 0.6955)$, also ordered according to the order in $c^T$.

The mean first passage time matrix is given by

$$M = \begin{bmatrix} 3.1097 & 15.2435 & 20.0601 & 20.1581 & 55.7987 \\ 9.5987 & 3.6974 & 22.3742 & 23.2789 & 57.7567 \\ 8.8444 & 17.0326 & 5.4278 & 22.1001 & 56.8645 \\ 7.6091 & 16.3412 & 21.0051 & 6.7752 & 56.5528 \\ 9.0204 & 17.6672 & 22.1062 & 22.2926 & 13.1345 \end{bmatrix}.$$

There is no ordered relationship within the row sums with $(m_{\bullet 1}, m_{\bullet 2}, m_{\bullet 3}, m_{\bullet 4}, m_{\bullet 5}) = (114.3702, 116.7059, 110.2695, 108.2834, 84.2210)$ but the vector of column sums is $(m_{1 \bullet}, m_{2 \bullet}, m_{3 \bullet}, m_{4 \bullet}, m_{5 \bullet}) = (38.1824, 69.9820, 90.9734, 94.6048, 240.1074)$, all with $c_i \geq c_j$ implying $m_{i \bullet} \leq m_{j \bullet}$ for $i \leq j$. We have not been able to establish any general results of such a nature for general finite Markov chains. Kemeny's constant for this chain is 16.042.

This example poses some unexpected inequalities that are not true in general, as seen in the following example.

**Example 4: (Eight state Markov chain)**

Funderlic and Meyer ([1]) provide an example involving the analysis of radiophosphorous kinetics in an aquarium system. This leads to a Markov chain with eight states. The states have been reordered so that the transition matrix has column sums with $c_i > c_j$ for each state $i < j$.

$$P = \begin{bmatrix} 0.478 & 0.270 & 0 & 0 & 0.150 & 0 & 0.055 & 0.047 \\ 0.130 & 0.870 & 0 & 0 & 0 & 0 & 0 & 0 \\ 0.320 & 0 & 0.669 & 0.011 & 0 & 0 & 0 & 0 \\ 0.088 & 0 & 0 & 0.912 & 0 & 0 & 0 & 0 \\ 0.150 & 0 & 0 & 0 & 0.740 & 0.110 & 0 & 0 \\ 0.300 & 0 & 0 & 0.011 & 0 & 0.689 & 0 & 0 \\ 0.260 & 0 & 0 & 0 & 0 & 0 & 0.740 & 0 \\ 0.600 & 0 & 0.400 & 0 & 0 & 0 & 0 & 0 \end{bmatrix},$$

with column sum vector $c^T = (2.326, 1.140, 1.069, 0.934, 0.890, 0.799, 0.795, 0.047)$, and stationary probability vector $\pi^T =$
$(0.2378, 0.4938, 0.0135, 0.0078, 0.1372, 0.0485, 0.0503, 0.0112)$.

In this example the ordering of the stationary probabilities does not follow that of the column sums with, for example, $\pi_1 < \pi_2$ even though $c_1 > c_2$.

The *H* matrix is given by

$$H = \begin{bmatrix} 1.0361 & 1.0557 & -0.3520 & -1.4033 & 0.1698 & -0.2611 & -0.1630 & 0.0428 \\ -0.7920 & 4.9494 & -0.4558 & -1.4630 & -0.8853 & -0.6343 & -0.5499 & -0.0431 \\ 0.2280 & -0.6227 & 2.6233 & -1.0520 & -0.2964 & -0.4260 & -0.3340 & 0.0048 \\ -1.6658 & -4.5559 & -0.5054 & 9.8722 & -1.3889 & -0.8124 & -0.7346 & -0.0842 \\ -0.2423 & -1.5994 & -0.4246 & -1.2750 & 3.2785 & 0.8385 & -0.4335 & -0.0173 \\ 0.1760 & -0.7306 & -0.4008 & -1.0294 & -0.3264 & 2.7789 & -0.3450 & 0.0024 \\ 0.1216 & -0.8436 & -0.4039 & -1.4332 & -0.3577 & -0.4477 & 3.4897 & -0.0002 \\ 0.4751 & -0.1095 & 0.8246 & -1.2706 & -0.1538 & -0.3755 & -0.2817 & 1.0165 \end{bmatrix}$$

The row sums of *H* are all the same at the value 0.125, ( =1/8, consistent with (2)). The column sums do not exhibit any pattern except through the relationship given by Theorem 2(b). The non-negativity of the diagonal elements $h_{jj}$ is consistent with (29) but the variation of the signs associated with the off-diagonal elements bears no discernible pattern.

The mean first passage time matrix *M* is given by

$$M = \begin{bmatrix} 4.21 & 7.88 & 220.32 & 1454.40 & 22.66 & 62.66 & 72.62 & 87.13 \\ 7.69 & 2.03 & 228.01 & 1462.09 & 30.35 & 70.35 & 80.32 & 94.82 \\ 3.40 & 11.28 & 74.05 & 1409.09 & 26.06 & 66.06 & 76.02 & 90.53 \\ 11.36 & 19.25 & 231.68 & 128.99 & 34.03 & 74.02 & 83.99 & 98.49 \\ 5.38 & 13.26 & 225.69 & 1437.84 & 7.29 & 39.99 & 78.00 & 92.50 \\ 3.62 & 11.50 & 223.93 & 1406.17 & 26.23 & 20.61 & 76.24 & 90.74 \\ 3.85 & 11.73 & 224.16 & 1458.24 & 26.51 & 66.50 & 19.88 & 90.97 \\ 2.36 & 10.24 & 133.19 & 1437.27 & 25.02 & 65.02 & 74.98 & 89.49 \end{bmatrix}.$$

As in the previous example there is no ordered relationship within the row sums of the mean first passage times, nor within the column sums. Kemeny's constant for this chain is 29.9194.

## 8. Summary

By introducing the matrix $H = [h_{ij}] = [I - P + ec^T]^{-1}$ where $c^T$ is the row vector of column sums of the transition matrix $P$, the stationary probability of being in state $j$ can be expressed in terms of the elements of $H$ and $c^T$ as $\pi_j = \sum_{i=1}^{m} c_i h_{ij}$ for all $j = 1, 2, ..., m$. It is also shown that $H$ can be expressed in terms of $Z$, Kemeny and Snell's fundamental matrix, thereby comparing the properties of the two matrices. The mean first passage times can also be expressed in the terms of the elements of $H$ and $c^T$ as $m_{ii} = 1 / \sum_{i=1}^{m} c_i h_{ij}$ and $m_{ij} = (h_{jj} - h_{ij}) / (\sum_{i=1}^{m} c_i h_{ij})$ for $i \neq j$. Some new relationship connecting the $\{c_i\}$ and $m_{ij}$ were derived; in particular $\pi_j = c_j / (m - \sum_{i=1}^{m} m_{ij} + \sum_{i=1}^{m} c_i m_{ij}) = 1 / (m - \sum_{i \neq j} m_{ij} + \sum_{i \neq j} c_i m_{ij})$. Inter-relationships between these aforementioned quantities were explored, including Kemenys' constant. Some general inequalities, based upon knowledge of the $\{c_i\}$ were explored and a variety of examples considered. Universal relationships were not achieved but some useful relationships were explored leaving some unanswered questions with scope for further investigations.

## References

[1] R. Funderlic, C.D. Meyer, Sensitivity of the stationary distribution vector for an ergodic Markov chain, Linear Algebra Appl. 76 (1986) 1-17.

[2] J. J. Hunter, Generalized inverses and their application to applied probability problems, Linear Algebra Appl. 45 (1982) 157-198.

[3] J. J. Hunter, Mathematical Techniques of Applied Probability, Volume 1, Discrete Time Models: Basic Theory, Academic, New York, 1983.

[4] J. J. Hunter, Mathematical Techniques of Applied Probability, Volume 2, Discrete Time Models: Techniques and Applications, Academic, New York, 1983.

[5] J. J. Hunter, Parametric Forms for Generalized inverses of Markovian Kernels and their Applications, Linear Algebra Appl. 127 (1990) 71-84.

[6] J. J. Hunter, Mixing times with applications to perturbed Markov chains, Linear Algebra Appl. 417 (2006) 108-123.

[7] J. J. Hunter, Simple procedures for finding mean first passage times in Markov chains, Asia-Pacific Journal of Operational Research, 24 (6), (2007), 813-829.

[8] J. J. Hunter, Variances of First Passage Times in a Markov chain with applications to Mixing Times, Linear Algebra Appl. 429 (2008) 1135-1162

[9] J. J. Hunter, Some stochastic properties of "semi-magic" and "magic" Markov chains, Linear Algebra Appl. 433 (2010), 893-907

[10] J. G. Kemeny, J.L. Snell, Finite Markov Chains, Van Nostrand, New York, 1960.

[11] S. Kirkland, Subdominant eigenvalues for stochastic matrices with given column sums, Electronic Journal of Linear Algebra, 18 (2009), 784-800. ISSN 1081-3810

[12] S. Kirkland, Column sums and the conditioning of stationary distribution for a stochastic matrix, Operators and Matrices, 4 (2010), 431-443. ISSN 1846-3886

[13] C. D. Meyer Jr., The role of the group generalized inverse in the theory of finite Markov chains, SIAM Rev. 17 (1975) 443-464.